\documentclass[12pt]{amsart}

\usepackage{amssymb,amsmath}
\usepackage[all]{xy}

\usepackage{amsmath}
\usepackage{amssymb}
\usepackage{url}
\usepackage{amscd}
\usepackage{mathrsfs}
\usepackage[dvips]{graphicx}

\usepackage[OT2,T1]{fontenc}
\DeclareSymbolFont{cyrletters}{OT2}{wncyr}{m}{n}
\DeclareMathSymbol{\Sha}{\mathalpha}{cyrletters}{"58}

\theoremstyle{plain}
\newtheorem{thm}{Theorem}

\newcommand{\Z}{{\mathbb Z}}

\newcommand{\Q}{{\mathbb Q}}
\newcommand{\R}{{\mathbb R}}

\newcommand{\Sela}{\mathrm{Sel}^{(\phi)}(E^{A,B}_r/\Q)}
\newcommand{\Selaa}{\mathrm{Sel}^{(\phi)}(E^{a,b}_r/\Q)}
\newcommand{\Selb}{\mathrm{Sel}^{(\widehat{\phi})}(\widehat{E}^{A,B}_r/\Q)}
\newcommand{\Selbb}{\mathrm{Sel}^{(\widehat{\phi})}(\widehat{E}^{a,b}_r/\Q)}

\begin{document}


\title[On rank zero twists of elliptic curves]
      {On positive proportion of rank zero twists of elliptic curves over $\Q$}



\author[Xiong]{Maosheng Xiong$^{\star}$}


\address{Maosheng Xiong: Department of Mathematics \\
Hong Kong University of Science and Technology \\
Clear Water Bay, Kowloon\\
P. R. China}
\email{mamsxiong@ust.hk}



\thanks{Maosheng Xiong is supported by RGC grant number 606211 and SBI12SC05 from Hong Kong.}

\begin{abstract}
Extending the idea of \cite{dab2} and using the 2-descent method, we provide three general families of elliptic curves over $\Q$ such that a positive proportion of prime-twists of such elliptic curves have rank zero simultaneously.

\end{abstract}

\maketitle

\section{Introduction}

Let $E$ be an elliptic curve over $\Q$ given by the Weierstrass equation $y^2=x^3+ax^2+bx+c \, (a,b,c \in \Z)$. For any square-free integer $d$, the $d$-th quadratic twist $E_d$ of $E$ is the elliptic curve given by the equation $y^2=x^3+adx^2+bd^2x+cd^3$. Denote by $r_{E_d}$ the rank of the Mordell-Weil group $E_d(\Q)$.

Statistics of rank zero quadratic twists of elliptic curves have been an interesting subject of study for some time. It is known from the work of Waldspurger \cite{wal} (combined with the work of Kolyvagin \cite{kol}, and of Wiles and others \cite{wil}) that $r_{E_d}=0$ for infinitely many square-free $d$'s. Hoffstein and Luo \cite{hof} proved that, for any fixed $E$, there exist infinitely many odd square-free $d$ with no greater than $3$ prime factors such that $r_{E_d}=0$. Ono and Skinner \cite{ono} proved that, for any fixed $E$, $|\{|d| \le X: r_{E_d}=0\}| \gg X/\log X$. On the other hand, it is believed \cite{gol} that a positive proportion of twists $E_d$ have rank zero, and this was proved by Iwaniec and Sarnak \cite{iwa2} under the Riemann hypothesis. Unconditionally, such positive proportion results are only known for a few specific curves \cite{jam, koh,vat,yu4}. Ono and Skinner \cite{ono} proved that, when $E$ has conductor $\le 100$, then $E_p$ or $E_{-p}$ has rank zero for a positive proportion of primes $p$.

In an interesting and beautiful paper \cite{dab2} Dabrowski proved the following result.
\begin{thm} \label{thm1} For any positive integer $k$ there are pairwise non-isogenous elliptic curves $E^1,\ldots,E^k$ such that $r_{{E^1}_p}=\cdots=r_{{E^k}_p}=0$ for a positive proportion of primes $p$.
\end{thm}

The ingenious idea of the proof of Theorem \ref{thm1} \cite{dab2} is based on the 2-descent method applied to the explicit family of elliptic curves $E^{A,-B}: y^2=x(x+A)(x-B)$ where $A$ is a prime, $B$ is a prime or a product of two primes such that $A+B=2^{2m}$ for some positive integer $m$. Chen's theorem \cite{chen} is required here to show that there are infinitely many such elliptic curves $E^{A,-B}$. As was remarked by Dabrowski \cite{dab2}, Theorem 1 can be also be proved by using the 2-descent method on the Setzer-Neumann curves, and the existence of infinitely many such elliptic curves was guaranteed by a result of Iwaniec \cite{iwa1}. Previously Dabrowski and Wieczorek \cite{dab1} proved Theorem \ref{thm1} by using the 2-descent method on the specific elliptic curves $y^2=x(x-2^m)(x+q-2^m)$ and by assuming the twin prime conjecture. 

The purpose of this paper is to show that there is an abundance of elliptic curves which could be used to prove Theorem \ref{thm1}.
Define
\begin{eqnarray} \label{1:eab} E^{A,B}: y^2=x(x+A)(x+B).\end{eqnarray}
We first prove

\begin{thm} \label{thm2} Let $|A|,|B|$ be primes such that
\begin{itemize}
\item[(1).] $A \equiv 1 \pmod{8}, B \equiv 3 \pmod{8}$,

\item[(2).] $A+B \ge 0$ or $AB<0$,
\end{itemize}
and let $r$ be any prime number coprime to $AB(A-B)$ and satisfy the following three conditions:
\begin{itemize}
\item[(i).] $r \equiv 7 \pmod{8}$,

\item[(ii).] $\left(\frac{A}{r}\right)=\left(\frac{B}{r}\right)=-1$,

\item[(iii).] For any odd prime $p|(A-B)$, we have $\left(\frac{-Br}{p}\right)=-1$.

\end{itemize}
Here $\left(\frac{\cdot}{\cdot}\right)$ is the Legendre symbol. Then $r_{E^{A,B}_r}=0$.
\end{thm}

Now Theorem \ref{thm1} can be proved easily: take distinct primes $p_1, \ldots,p_r$ and $q$ with $p_i \equiv 1 \pmod{8} \forall \, i$ and $q \equiv 5 \pmod{8}$, and consider the $r$-twist of elliptic curves $E^{p_i,q}$. Let $r$ be a prime number coprime to $p_iq(p_i-q) \, \forall i$ and satisfy
\begin{itemize}
\item[(i).] $r \equiv 7 \pmod{8}$,

\item[(ii).] $\left(\frac{q}{r}\right)=-1$ and $\left(\frac{p_i}{r}\right)=-1 \, \forall \, i$,

\item[(iii).] For any odd prime $p|\prod_i(p_i-q)$, we have $\left(\frac{-qr}{p}\right)=-1$.

\end{itemize}
Then $r_{E^{p_i,q}_r}=0$ for all $i$. Chinese Remainder Theorem and Dirichlet's Theorem on primes in arithmetic progressions clearly show a positive proportion of primes $r$ satisfying the above conditions.

The proof of Theorem \ref{thm2} is also based on the 2-descent method. Consider the two-isogeny $\phi: E^{A,B}_r \to \widehat{E}^{A,B}_r:Y^2=X^3-2(A+B)rX^2+(A-B)^2r^2X$, defined by $\phi((x,y))=(y^2/x^2,-y(x^2-ABr^2)/x^2)$, and let $\widehat{\phi}$ denote the dual isogeny. Under the assumptions of Theorem \ref{thm2}, we actually prove that the Selmer groups satisfy $\mathrm{Sel}^{(\phi)}(E^{A,B}_r/\Q) \simeq \{0\}$ and $\mathrm{Sel}^{(\widehat{\phi})}(\widehat{E}^{A,B}_r/\Q) \simeq (\Z/2\Z)^2$. Then Theorem \ref{thm2} follows from the fundamental formula (\cite[p. 314]{sil})
\begin{eqnarray*} r_{E^{A,B}_r} &=& \dim_2\mathrm{Sel}^{(\phi)}(E^{A,B}_r/\Q)+\dim_2 \mathrm{Sel}^{(\widehat{\phi})}(\widehat{E}^{A,B}_r/\Q)\\
&& -
\dim_2\Sha(E^{A,B}_r/\Q)[\phi]-\dim_2 \Sha(\widehat{E}^{A,B}_r/\Q)[\widehat{\phi}]-2.  \end{eqnarray*}

Theorem \ref{thm2} is about rank zero $r$-twists of $E^{A,B}$ for $|A|$ and $|B|$ being primes. If $A,B$ are some integers in general, it is still possible to find rank zero $r$-twists of $E^{A,B}$, given the elementary nature of the 2-descent method, however, the conditions on such $r$'s are too complicated to write down. On the other hand, we prove

\begin{thm} \label{thm3} Let $A,B$ be integers such that $A \equiv 1 \pmod{8}, B \equiv 3 \pmod{8}$, and for $E^{A,B}$, we assume that
\begin{eqnarray} \label{1:sel0} \mathrm{Sel}^{(\phi)}(E^{A,B}/\Q) \simeq \{0\}, \quad \mathrm{Sel}^{(\widehat{\phi})}(\widehat{E}^{A,B}/\Q) \simeq (\Z/2\Z)^2.\end{eqnarray}
Then for any prime number $r$ coprime to $AB(A-B)$ such that
\begin{itemize}
\item[(i).] $r \equiv 7 \pmod{8}$, and

\item[(ii).] for any odd prime $p|AB(A-B)$, we have $\left(\frac{r}{p}\right)=1$,
\end{itemize}
we have $r_{E^{A,B}_r}=0$.
\end{thm}

We remark that Condition (\ref{1:sel0}) is natural and easy to check for $E^{A,B}$. It also implies that $r_{E^{A,B}}=0$. Many families of rank zero elliptic curves have been found by this way (see for example \cite{fxi,dab1,dab2,fxu}). Moreover, via infinitely many elliptic curves $E^{A,B}$ satisfying Condition (\ref{1:sel0}), Theorem \ref{thm1} can be proved easily.

The elliptic curve $E^{A,B}$ given in (\ref{1:eab}) can be characterized as having full 2-torsion points over $\Q$, that is, $E^{A,B}(\Q)[2] \simeq (\Z/2\Z)^2$. Now we consider elliptic curves with one non-trivial 2-torsion point over $\Q$, that is, elliptic curves $E^{a,b}$ ($a,b \in \Z$) given by the equation
\begin{eqnarray} \label{1:eab2}
E^{a,b}: y^2=x(x^2+ax+b).
\end{eqnarray}
Consider the two-isogeny $\phi: E^{a,b}_r \to \widehat{E}^{a,b}_r:Y^2=X^3-2arX^2+(a^2-4b)r^2X$, defined by $\phi((x,y))=(y^2/x^2,-y(x^2-br^2)/x^2)$, and let $\widehat{\phi}$ denote the dual isogeny. We also have the fundamental formula (\cite[p. 314]{sil})
\begin{eqnarray*} r_{E^{a,b}_r} &=& \dim_2\mathrm{Sel}^{(\phi)}(E^{a,b}_r/\Q)+\dim_2 \mathrm{Sel}^{(\widehat{\phi})}(\widehat{E}^{a,b}_r/\Q) \\
&& -
\dim_2\Sha(E^{a,b}_r/\Q)[\phi]-\dim_2 \Sha(\widehat{E}^{a,b}_r/\Q)[\widehat{\phi}]-2.  \end{eqnarray*}
We assume that neither $a^2-4b$ nor $b$ is a perfect square, because otherwise, either $E^{a,b}$ or $\widehat{E}^{a,b}$ would have full 2-torsion points over $\Q$, hence reducing to elliptic curves we have considered before. We prove

\begin{thm} \label{thm4} Let $a,b \in \Z$ such that neither $a^2-4b$ nor $b$ is a perfect square. For $E^{a,b}$ given in (\ref{1:eab2}), we assume that
\begin{eqnarray} \label{1:sel2} \mathrm{Sel}^{(\phi)}(E^{a,b}/\Q) \simeq \Z/2\Z, \quad \mathrm{Sel}^{(\widehat{\phi})}(\widehat{E}^{a,b}/\Q) \simeq \Z/2\Z.\end{eqnarray}
Then for any prime number $r$ coprime to $b(a^2-4b)$ such that
\begin{itemize}
\item[(i).] $r \equiv 7 \pmod{8}$,

\item[(ii).] for any odd prime $p|b(a^2-4b)$, we have $\left(\frac{r}{p}\right)=1$,

\item[(iii).] either
\begin{eqnarray} \label{1:leg} \left(\frac{b}{r}\right)= \left(\frac{a^2-4b}{r}\right)=-1,\end{eqnarray}
or
\[\left(\frac{b}{r}\right)= \left(\frac{a^2-4b}{r}\right)=-\left(\frac{a+\sqrt{b}}{r}\right)=
    -\left(\frac{a+\sqrt{a^2-4b}}{r}\right)=1,\]
\end{itemize}
we have $r_{E^{a,b}_r}=0$.
\end{thm}
We remark that if $a=2a', b=2b'$ and $a' \equiv b' \equiv 3 \pmod{4}$, then Condition (\ref{1:leg}) is satisfied. On the other hand, Condition (\ref{1:sel2}) is natural and easy to check for $E^{a,b}$ as well, and many rank zero elliptic curves have been found in this way (\cite{fxi2}). Theorem \ref{thm1} can also be proved via infinitely many elliptic curves $E^{a,b}$ satisfying the conditions of Theorem \ref{thm4}.

The main ingredient of the proofs of Theorems \ref{thm2}--\ref{thm4} is the 2-descent method applied to the elliptic curves considered above. We prove Theorems \ref{thm2} and \ref{thm3} in Section 3, and prove Theorem \ref{thm4} in Sections 5.


\section{Proof of Theorems \ref{thm2} and \ref{thm3}}

The 2-descent method is explained in the last chapter of Silverman's book (\cite{sil}, see also \cite{ata,dab2}). For clarity we specify the 2-descent method for elliptic curves $E^{A,B}$ given by (\ref{1:eab}) below.

\subsection{2-descent and $E^{A,B}$}

For any integer $M$, let $\Sigma(M)$ be the set of prime numbers dividing $M$, and let $\triangle(M)$ be set of (positive or negative) square-free divisors of $M$. Let
\[
C^{(r)}_{d}: dw^2=t^4-2(A+B)\frac{r}{d}t^2z^2+(A-B)^2\frac{r^2}{d^2}z^4\,
, \]
\[
C^{'(r)}_{d}: dw^2=t^4+(A+B)\frac{r}{d}t^2z^2+AB \frac{r^2}{d^2} z^4\,
\]
be the principal homogeneous spaces under the actions of the elliptic curves $E^{A,B}_r$ and $\widehat{E}^{A,B}_r$ previously defined. Using \cite[Proposition 4.9, p. 302]{sil}, we have the following identifications:
\[ \Sela \simeq \left\{ d \in \triangle((A-B)r): C^{(r)}_{d}(\Q_v)
\ne \emptyset \, \forall v \in \Sigma(2AB(A-B)r) \cup \{\infty\} \right\}\, , \]
\[ \Selb \simeq \left\{ d \in \triangle(ABr): C^{'(r)}_{d}(\Q_v)
\ne \emptyset \, \forall v \in \Sigma(2AB(A-B)r) \cup \{\infty\} \right\} \, , \]
where $C^{(r)}_{d}(\Q_v)
\ne \emptyset$ (or $C^{'(r)}_{d}(\Q_v) \ne \emptyset$) means that $C^{(r)}_{d}$ (or $C^{'(r)}_{d}$) has non-trivial solutions
$(w,t,z) \ne (0,0,0)$ in $\Q_v$. We know that
\[ \{1 \} \subseteq \Sela, \qquad  \{1,AB,-Ar,-Br \} \subseteq \Selb,  \]
since each of the corresponding homogeneous spaces has
non-trivial solutions in $\Q$.

\noindent {\bf Proof of Theorem \ref{thm2}.} Let $|A|,|B|$ and $r$ be primes under the conditions of Theorem \ref{thm2}. We first prove that $\Sela = \{1\}$. For any $d \in \triangle((A-B)r)$, if $d<0$, since $A+B \ge 0$ or $AB<0$, clearly $C^{(r)}_d(\R) =\emptyset$. If $2|d$, let $(w,t,z) \ne (0,0,0)$ be a solution of $C^{(r)}_d$ in $\Q_2$. Since $A \equiv 1 \pmod{8}, B \equiv 5 \pmod{8}$, considering $C^{(r)}_d$, at least two numbers in the set $\{1+2v_2(w),4v_2(t),1+2v_2(t)+2v_2(z),2+4v_2(z) \}$ reach the same minimal value, here for any prime $p$, we denote by $v_p$ the standard $p$-adic exponential valuation. Clearly these two numbers must be $1+2v_2(w)$ and $1+2v_2(t)+2v_2(z)$, which is impossible. Now we need to consider $d > 1$ such that there is an odd prime $p|(A-B)$ with $p|d$. Let $(w,t,z) \ne (0,0,0)$ be a solution of $C^{(r)}_d$ in $\Q_p$. Then at least two numbers in the set $\{1+2v_p(w),4v_p(t),-1+2v_p(t)+2v_p(z),2v_p(A-B)-2+4v_p(z) \}$ reach the same minimal value. These two numbers must be $1+2v_p(w)$ and $-1+2v_p(t)+2v_p(z)$. Hence the equation
\[\frac{d}{p}w^2 \equiv -2(A+B) \frac{r}{d/p}t^2z^2 \pmod{p}\]
must be solvable in ${\Z_p}^*:=\Z_p-p\Z_p$, where $\Z_p$ is the set of $p$-adic integers. By Hensel's lemma, this implies that
\[1=\left(\frac{-2(A+B)r}{p}\right)=\left(\frac{-2(A-B+2B)r}{p}\right)=
\left(\frac{-rB}{p}\right),\]
which contradicts Condition (iii) of Theorem \ref{thm2}.

Finally if $d=r$, let $(w,t,z) \ne (0,0,0)$ be a solution of $C^{(r)}_d$ in $\Q_r$. Then at least two numbers in the set $\{1+2v_r(w),4v_r(t),2v_r(t)+2v_r(z),4v_r(z) \}$ reach the same minimal value. We may assume that $v_r(t)=v_r(z)=0$ and $v_r(w) \ge 0$, hence
\[t^4-2(A+B)t^2z^2+(A-B)^2z^4 \equiv 0 \pmod{r}\]
is solvable in ${\Z_r}^*$. That is,
\[(t^2-(A+B)z^2)^2 \equiv 4AB z^4 \pmod{r}. \]
Since $\left(\frac{AB}{r}\right)=1$, we have
\[t^2-(A+B-2\sqrt{AB})z^2 \equiv 0 \pmod{r}, \]
or
\[t^2-(A+B+2\sqrt{AB})z^2 \equiv 0 \pmod{r}. \]
This is not possible by Hensel's lemma and by Conditions (i)(ii) of Theorem \ref{thm2}, since 
\[\left(\frac{A+B \pm 2\sqrt{AB}r}{r}\right)=\left(\frac{-1}{r}\right)\left(\frac{(\sqrt{-A} \pm \sqrt{-B})^2}{r}\right)=-1.\]
This shows that $\Sela = \{1\}$.

Next we prove that $\Selb=\{1,AB,-Ar,-Br\}$. For $d \in \triangle(ABr)$, for $d=-1$, we have
\[C^{'(r)}_{-1}: \, -w^2=t^4-(A+B)rt^2z^2+ABr^2z^4. \]
Since $\left(\frac{-1}{r}\right)=-1$ and $\left(\frac{-AB}{r}\right)=-1$, by Hensel's Lemma, we find that $C^{'(r)}_{-1}(\Q_r)=\emptyset$. For $d=A$, we have
\[C^{'(r)}_{A}: \, Aw^2=t^4+(A+B)\frac{r}{A}t^2z^2+\frac{Br^2}{A}z^4. \]
Since $\left(\frac{A}{r}\right)=\left(\frac{B}{r}\right)=-1$, by Hensel's Lemma, we find that $C^{'(r)}_{A}(\Q_r)=\emptyset$. For $d=r$, we have
\[C^{'(r)}_{r}: \, rw^2=t^4+(A+B)t^2z^2+ABz^4. \]
Solving it in $\Q_2$, at least two numbers in the set $\{2v_2(w),4v_2(t),1+2v_2(t)+2v_2(z),4v_2(z) \}$ reach the same minimal value, which may be assumed to be zero. This implies that at least one of the equations
\[rw^2 \equiv t^4 \pmod{8}, \quad rw^2 \equiv ABz^4 \pmod{8}, \]
and
\[rw^2 \equiv t^4+(A+B)t^2z^2+ ABz^4=(t^2+Az^2)(t^2+Bz^2) \pmod{8} \]
is solvable in ${\Z_2}^*$. This is impossible by Conditions (1) and (i) of Theorem \ref{thm2}. Since $\Selb \subset \triangle(ABr)$ is a group, we conclude that $\Selb =\{1,AB,-Ar,-Br\}$. This shows $r_{E^{A,B}_r}=0$, completing the proof of Theorem \ref{thm2}. \quad $\square$

\noindent {\bf Proof of Theorem \ref{thm3}.} For any $d' \in \triangle(A-B)$, let $d=d'$ or $d=d'r$. If $d \in \Sela$, then $C^{(r)}_d(\Q_v) \ne \emptyset$ for any $v \in \Sigma(2AB(A-B)) \cup \{\infty\}$. Since by (i) and (ii) of Theorem \ref{thm3}, $r$ is a square in any of such $\Q_v$'s, we find $d' \in \mathrm{Sel}^{\phi}(E^{A,B}/\Q)$. Hence $d'=1$ from Condition (2) of Theorem \ref{thm3}. For $d=r$, since $A B \equiv 3 \pmod{4}$ and by (ii) of Theorem \ref{thm3} we have
\[\left(\frac{AB}{r}\right)=-1. \]
From the proof of Theorem \ref{thm2}, we find $C^{(r)}_r(\Q_r) = \emptyset$. Hence $\Sela=\{1\}$.

For any $d' \in \triangle(AB)$, let $d=d'$ or $d=d'r$. If $d \in \Selb$, since $r$ is a square in $\Q_v$ for any $v \in \triangle(2AB(A-B))$, by Condition (2) of Theorem \ref{thm3} we find $d' \in \{1,AB,-A,-B\}$. Hence $\Selb \subseteq \{1,AB,-A,-B,r,ABr,-Ar,-Br\}$. 

For $d=r$, from the proof of Theorem \ref{thm2}, we find $C^{'(r)}_{r}(\Q_2) = \emptyset$. Using that $\Selb$ is a group, we conclude that $\Selb =\{1,AB,-Ar,-Br\}$. This shows $r_{E^{A,B}_r}=0$, which completes the proof of Theorem \ref{thm3}. \quad $\square$

\section{Proof of Theorem \ref{thm4}}
For clarity we specify the 2-descent method for elliptic curves $E^{a,b}$ given by (\ref{1:eab2}) below (see \cite{sil,ata}).
\subsection{2-descent and $E^{a,b}$}

Let
\[
C^{(r)}_{d}: dw^2=t^4-2a\frac{r}{d}t^2z^2+(a^2-4b)\frac{r^2}{d^2}z^4\,
, \]
\[
C^{'(r)}_{d}: dw^2=t^4+\frac{ar}{d}t^2z^2+b \frac{r^2}{d^2} z^4\,
\]
be the principal homogeneous spaces under the actions of the elliptic curves $E^{a,b}_r$ and $\widehat{E}^{a,b}_r$ previously defined. Using \cite[Proposition 4.9, p. 302]{sil}, we have the following identifications:
\[ \Selaa \simeq \left\{ d \in \triangle((a^2-4b)r): C^{(r)}_{d}(\Q_v)
\ne \emptyset \, \forall v \in \Sigma(2b(a^2-4b)r) \cup \{\infty\} \right\}\, , \]
\[ \Selbb \simeq \left\{ d \in \triangle(br): C^{'(r)}_{d}(\Q_v)
\ne \emptyset \, \forall v \in \Sigma(2b(a^2-4b)r) \cup \{\infty\} \right\} \, , \]
where $C^{(r)}_{d}(\Q_v)
\ne \emptyset$ ( or $C^{'(r)}_{d}(\Q_v) \ne \emptyset$) means that $C^{(r)}_{d}$ (or $C^{'(r)}_{d}$) has non-trivial solutions
$(w,t,z) \ne (0,0,0)$ in $\Q_v$. We know that
\[ \{1, a^2-4b\} \subseteq \Selaa, \qquad  \{1,b \} \subseteq \Selbb,  \]
since each of the corresponding homogeneous spaces has
non-trivial solutions in $\Q$.

\subsection{Proof of Theorem \ref{thm4}} Similar to the proof of Theorem \ref{thm3}, since $r$ is a square in $\Q_v$ for any $v \in \triangle(2b(a^2-4b))$, by Condition (4) of Theorem \ref{thm4}, we have
\[ \Selaa \subseteq \{1,a^2-4b,r,(a^2-4b)r \}, \qquad  \Selbb \subseteq \{1,b,r,br\}.\]
It suffices to prove that
\[ r \notin \Selaa, \qquad  r \notin \Selbb.\]
For
\[
C^{(r)}_{r}: rw^2=t^4-2at^2z^2+(a^2-4b)z^4\, , \]
let $(w,t,z) \ne (0,0,0)$ be a solution of $C^{(r)}_r$ in $\Q_r$. Then at least two numbers in the set $\{1+2v_r(w),4v_r(t),v_r(a)+2v_r(t)+2v_r(z),4v_r(z)\}$ reach the minimal value, which we may assume to be zero. Hence $v_r(t)=v_r(z)=0,v_r(w) \ge 0$. So
\[(t^2-az^2)^2 \equiv 4b z^4 \pmod{r}\]
is solvable in ${\Z_r}^*$. This requires that $\left(\frac{b}{r}\right)=1$. Moreover, it implies that at least one of the equations
\[t^2 \equiv (a+2 \sqrt{b})z^2 \pmod{r},\]
\[t^2 \equiv (a-2 \sqrt{b})z^2 \pmod{r},\]
is solvable in ${\Z_r}^*$. This means by Hensel's lemma
\[\left(\frac{a+2\sqrt{b}}{r}\right)=1, \quad \mbox{ or } \left(\frac{a-2\sqrt{b}}{r}\right)=1. \]
It is easy to see that this contradicts (iii) of Theorem \ref{thm4}. 

For
\[
C^{'(r)}_{r}: rw^2=t^4+at^2z^2+bz^4\, , \]
similar to the argument before, we see that $C^{'(r)}_{r}(\Q_r) = \emptyset$. This concludes the proof of Theorem \ref{thm4}. \quad $\square$

\end{document}